\documentclass[11pt, reqno]{amsart}     
\usepackage[dvipdfmx]{graphicx}
\usepackage{graphicx}
\usepackage{amssymb,amstext,amsfonts,amscd}
\usepackage{amsmath}
\usepackage{latexsym}
\usepackage{amsthm}
\usepackage{multicol}
\usepackage{fancybox}
\usepackage{here}
\usepackage{multirow}
\usepackage{mathrsfs}

\usepackage{enumitem}

\makeatletter
\@namedef{subjclassname@2020}{%
  \textup{2020} Mathematics Subject Classification}
\makeatother

\def\opn#1#2{\def#1{\operatorname{#2}}}
\opn\max{max}
\opn\min{min}
\opn\rank{rank}
\opn\Ker{Ker}
\opn\id{id}
\opn\mod{mod}
\opn\det{det}
\opn\Cone{Cone}
\opn\Int{Int}
\opn\sign{sign}

\newcommand{\bz}{\mathbb{Z}}

\newcommand{\br}{\mathbb{R}}
\newcommand{\bc}{\mathbb{C}}

\newcommand{\zero}{\boldsymbol{0}}
\newcommand{\ba}{\boldsymbol{a}}   

\newcommand{\bbu}{\boldsymbol{u}}

\newcommand{\bw}{\boldsymbol{w}}

\newcommand{\bbz}{\boldsymbol{z}}

\newcommand\De{\Delta}

\newcommand\rdeg{{\rm{rdeg}\/}}
\newcommand\pdeg{{\rm{pdeg}\/}}

\newcommand{\QED}{$\Box$}

\newtheorem{theorem}{Theorem}
\newtheorem{definition}[theorem]{Definition}

\newtheorem{proposition}[theorem]{Proposition}
\newtheorem{lemma}[theorem]{Lemma}

\newtheorem{example}[theorem]{Example}
\newtheorem{problem}[theorem]{Problem}

\numberwithin{equation}{section}

\begin{document}
\pagestyle{plain}

\title
{Toric resolutions of strongly mixed weighted homogeneous polynomial germs of type $J_{10}^-$}

\author{Sachiko Saito}

\dedicatory{In Commemoration of Professor Goo Ishikawa's all years of hard work}

\address{
Department of Mathematics Education, Asahikawa Campus, Hokkaido University of Education, 
Asahikawa 070-8621, Hokkaido, Japan}

\email{saito.sachiko@a.hokkyodai.ac.jp}



\begin{abstract}
We consider toric resolutions of 
some strongly mixed weighted homogeneous polynomials of type $J_{10}^-$. 
We show that the strongly mixed weighted homogeneous polynomial $f := f_{2,2,1,2,1,4}\ (k=3)$ (see \S \ref{section-J10}) has 
no mixed critical points on ${\bc^*}^2$ (Lemma \ref{no-mixed-critical-points}), and moreover, 
show that 
the strict transform $\tilde V$ of the mixed hypersurface singularity $V := f^{-1}(0)$ via the toric modification $\hat{\pi} : X \to \bc^2$, 
where we set $f := f_{2,2,1,2,1,4}\ (k=3)$, 
is not only a real analytic manifold outside of ${\tilde V} \cap \hat{\pi}^{-1}(\boldsymbol{0})$ 
but also a real analytic manifold as a germ of ${\tilde V} \cap \hat{\pi}^{-1}(\boldsymbol{0})$ (Theorem \ref{real-analytic-manifold}). 
\end{abstract}

\maketitle


\section{Introduction}

The famous Hilbert's 16th problem (see \cite{Gudkov74}, \cite{Wilson}, \cite{Ishi-Saito-Fukui2001} for example) asks the topology of 
nonsingular real algebraic curves of degree $6$ on the real projective plane $\br P^2$ and 
nonsingular real algebraic surfaces of degree $4$ in the real projective space $\br P^3$. 
D.A. Gudkov completed the isotopic classification (the 56 isotopy types) 
of nonsingular real algebraic curves of degree $6$ on $\br P^2$ in 1971 (cf. \cite{Gudkov74}). 
Especially he proved that the number of the isotopy types of nonsingular real algebraic $M$-curves of degree $6$ on $\br P^2$ 
is three. See Figure \ref{M-curves-of-degree} below: 
\begin{figure}[H]
\begin{center}
\includegraphics[width=10cm]{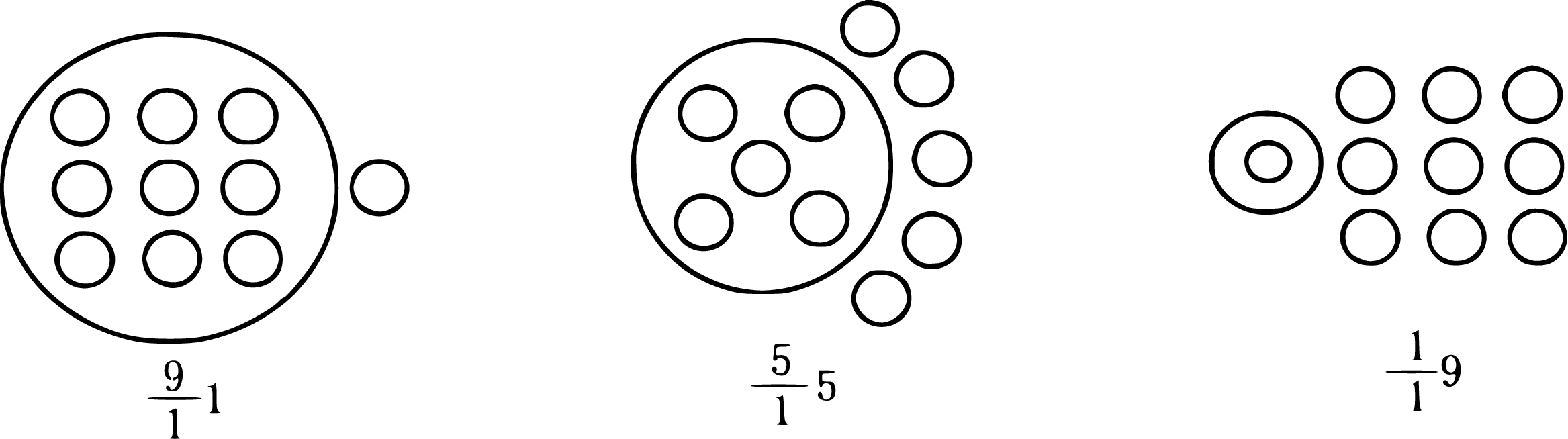}
\end{center}
\caption{The isotopy types of nonsingular real algebraic $M$-curves of degree $6$ on $\br P^2$}
\label{M-curves-of-degree}
\end{figure}
In the early 1970s, V. A. Rokhlin (\cite{Rokhlin72a}, \cite{Rokhlin72b}, \cite{Rokhlin73}) reproved this Gudkov's result by differential topology, 
which includes Smith theory on group actions, and some index theorems on $4$-dimensional manifolds. 
From complex geometric point of view, 
double coverings of $\bc P^2$ branched along nonsingular algebraic sextic curves 
and 
nonsingular algebraic quartic surfaces in $\bc P^3$ 
are {\em K3 surfaces}. 
V. M. Kharlamov (\cite{Kharlamov1976}) started applying the K3 surface theory and the deformation theory to Hilbert's 16th problem 
and obtained many remarkable results. 
Subsequently to Kharlamov's approaches, 
V.V. Nikulin (\cite{Nikulin79}) reproved Gudkov's isotopic classification of nonsingular real sextic curves 
by using the {\em moduli spaces} ({\em period domains}, or {\em Hilbert schemes}) of real projective K3 surfaces 
and his lattice (integral symmetric bilinear forms) theory in 1979. 

In the 1980s, O. Ya. Viro (\cite{Viro1984}) reconstructed Gudkov's 56 isotopy types of nonsingular real algebraic curves of degree $6$ 
by his {\em patchworking method} (see \cite{Viro2006} for the details), 
where he used the real weighted homogeneous polynomials of type $J_{10}^-$ 
with two variables: 
$$
f(z_1,z_2)  :=  (z_2 - z_1^2)(z_2 - 2z_1^2)(z_2 - k z_1^2), 
$$
where $k$ is a real number and $k>2$, and their nonsingular perturbations. 
He got all the nonsingular perturbations of the weighted homogeneous polynomials of type $J_{10}^-$ (\cite{Viro1989}). 
Roughly speaking, he patchworked two appropriate nonsingular perturbations of some polynomials of type $J_{10}^-$, 
and projectified the patchworked polynomial of degree $6$. 
Nonsingular perturbations whose zero sets (in $\br^2$) have many connected components 
are very useful for constructions of nonsingular real algebraic $M$-curves of degree $6$ on $\br P^2$. 
Thus, the weighted homogeneous polynomials of type $J_{10}^-$ were found to be very important for 
the constructions of nonsingular real algebraic curves of degree $6$ on $\br P^2$. 

In this paper 
we consider {\em mixed weighted homogeneous polynomials} of type $J_{10}^-$, 
especially, strongly mixed weighted homogeneous polynomials of type $J_{10}^-$. 
It would be expected that there exist some relations between 
the topology of the real parts of (toric) resolutions of mixed weighted homogeneous polynomials with real coefficients 
and 
that of nonsingular perturbations of real weighted homogeneous polynomials. 
In this paper 
we eventually show that 
the strongly mixed weighted homogeneous polynomial (the polynomial \eqref{IV-alpha3} in \S \ref{section-J10}):
$$
f := f_{2,2,1,2,1,4}\ (k=3)
$$
has no mixed critical points on ${\bc^*}^2$ (Lemma \ref{no-mixed-critical-points}), and moreover, 
show that 
the strict transform $\tilde V$ of the mixed hypersurface singularity $V := f^{-1}(0)$ via the toric modification $\hat{\pi} : X \to \bc^2$, 
where we set $f := f_{2,2,1,2,1,4}\ (k=3)$, 
is not only a real analytic manifold outside of $\tilde V \cap \hat{\pi}^{-1}(\zero)$ 
but also a real analytic manifold as a germ of $\tilde V \cap \hat{\pi}^{-1}(\zero)$ (Theorem \ref{real-analytic-manifold}). 

\bigskip

\section{Basic terminology about mixed polynomials}

Here we recall basic terminology about mixed polynomials. 

\subsection{Mixed functions and their radial Newton polyhedrons}\label{section-mixed-function}
Let $U$ be a neighborhood of $\zero$ in $\bc^n$ with $\bar{U}=U$, 
where $\bar{\bbz}$ stands for the complex conjugate $(\bar{z_1}, \dots , \bar{z_n})$ of $\bbz=(z_1, \dots , z_n) \in \bc^n$. 
Let $F(\bbz, \bw)$ be a complex valued holomorphic function on $U\times U$ with complex $2n$ variables. 
We assume that $F(\zero, \zero) = 0$
We define $f : U \to \bc$ by 
$$
f(\bbz, \bar{\bbz}) := F(\bbz, \bar{\bbz}). 
$$
We call such a $f$ (or $f(\bbz, \bar{\bbz})$) 
a {\em mixed analytic function} (or {\em mixed function}) on $U$. 
Let 
$
F(\bbz, \bw) = \sum_{\nu, \mu}c_{\nu, \mu} \bbz^\nu \bw^\mu
$
be the Taylor expansion of $F$ at $(\zero , \zero)$, 
where 
$
\nu = (\nu_1, \dots , \nu_n),\ \mu = (\mu_1, \dots , \mu_n),\ 
\nu_i \geq 0,\ \mu_j \geq 0,\ \ 
\bbz^\nu := z_1^{\nu_1}\cdots z_n^{\nu_n},\ \ \bw^\mu := w_1^{\mu_1}\cdots w_n^{\mu_n} .
$
Then we have 
\begin{equation}\label{Taylor}
f(\bbz, \bar{\bbz}) = \sum_{\nu, \mu}c_{\nu, \mu} \bbz^\nu \bar{\bbz}^\mu .
\end{equation}
Note that the coefficients of the Taylor expansion \eqref{Taylor} of a mixed function $f$ are unique. 
We call $f(\bbz, \bar{\bbz})$ a {\em mixed polynomial} 
if 
the number of monomials $c_{\nu, \mu} \bbz^\nu \bar{\bbz}^\mu,\ c_{\nu, \mu}\neq 0$ is finite. 

\begin{definition}
{\rm 
Let $f(\bbz, \bar{\bbz})$ be a mixed function on $U \ (\subset \bc^n)$. 
We say 
$
\ba = (a_1, \dots , a_n) \in U
$
is a {\em mixed critical point} (or a {\em mixed singular point}) of $f$ 
if 
the rank of the differential map 
$$
(df)_{\ba} : T_{\ba} \bc^n \to T_{f(\ba)} \bc \cong T_{f(\ba)} \br^2
$$
is less than $2$. 
We say $\ba \in U$ is a {\em mixed regular point} of $f$ 
if it is not a mixed critical point of $f$. 
}
\end{definition}

\medskip

We set 
$K_+^n := \{ (x_1, \dots , x_n) \in K^n \ |\ x_i \geq 0\ \text{for\ every}\ i \}$, where $K=\br \ \text{or} \ \bz$. 

For the germ $(f,\zero)$ (at $\zero \in \bc^n$) of a mixed function 
$
f(\bbz, \bar{\bbz}) = \displaystyle \sum_{\nu, \mu}c_{\nu, \mu} \bbz^\nu \bar{\bbz}^\mu
$, 
$$
\Gamma_+(f)
$$
is defined to be the convex hull of the set 
$
\displaystyle \bigcup_{c_{\nu, \mu} \neq 0}  (\nu + \mu) + \br_+^n
$. 
We call $\Gamma_+(f)$ 
the {\em (radial) Newton polyhedron} of $(f,\zero)$. 

For a ``weight vector" 
$
P={}^t(p_1, \dots , p_n)\ (\neq \zero) \in \bz_+^n
$, 
we define 
$d(P)$ 
to be the minimum value of the linear function 
$
P : \Gamma_+(f) \to \br,\ P(\xi) := \displaystyle \sum_{j=1}^n p_j \xi_j,
$
where 
$\xi = (\xi_1 , \dots , \xi_n) \in \Gamma_+(f)$. 
We set 
$$
\Delta(P) := \{ \xi \in \Gamma_+(f) \ |\ P(\xi) = d(P) \}, 
$$
which we call a {\em face} of $\Gamma_+(f)$. 
Note that $\Delta(P) \neq \emptyset$ by its definition. 

We say a weight vector $P={}^t(p_1, \dots , p_n) \in \bz_+^n$ is {\em strictly positive} ($P\gg 0$) if 
$p_i >0$ for every $i \ (=1, \dots , n)$. 
Note that a face $\Delta$ of $\Gamma_+(f)$ is compact 
if and only if 
$\Delta = \Delta(P)$ 
for some strictly positive weight vector $P$.    
For a \underline{compact} face $\Delta(P)$, we define 
$$
f_{\Delta(P)}(\bbz) \ (\text{or}\ f_P(\bbz)\,) := \sum_{\mu + \nu \in \Delta(P)}c_{\nu, \mu} \bbz^\nu \bar{\bbz}^\mu ,
$$
which we call a {\em face function {\rm (or}\ face polynomial{\rm )}} 
of a mixed function germ $(f,\zero)$ (\cite{Oka2018}, p.78).

\begin{definition}[\cite{Oka2018}, p.79]\label{convenient-function} 
{\rm 
Let $(f,\zero)$ be the germ of a mixed function $f$ at $\zero \in \bc^n$. 
For $I \subset \{1, \dots, n\}$, 
$f^I$ denotes the restriction of $f$ on the coordinate subspace $\bc^I := \{ \bbz \,|\, z_j=0,\, j\notin I \}$. 
A mixed function germ $(f,\zero)$ is called {\em convenient} 
if 
$f^{I} \not\equiv 0$ for every $I \subset \{ 1,2,\cdots,n \}$ with $I \neq \emptyset$. 
}
\end{definition}

\medskip

\subsection{Radially and polar weighted homogeneous polynomials}\label{section-whp}

\begin{definition}[\cite{Oka2018}, p.182;\ see also \cite{Oka2015}]
{\rm \ \ 

\begin{itemize}
\item A mixed polynomial 
$f(\bbz, \bar{\bbz}) = \sum_{\nu, \mu}c_{\nu, \mu} \bbz^\nu \bar{\bbz}^\mu$ is called 
{\em radially weighted homogeneous} if there exists a weight vector 
$P={}^t(p_1, \dots , p_n)\ (\neq \zero) \in \bz_+^n$   
and 
a positive integer $d_r \ (> 0)$ such that 
$$
c_{\nu, \mu} \neq 0 \Longrightarrow P(\nu + \mu) = \sum_{i=1}^n p_i(\nu_i + \mu_i) = d_r .
$$
We call $d_r$ {\em the radial degree} of $f$, and define 
$
\rdeg_P f := d_r .
$

\item A mixed polynomial $f(\bbz, \bar{\bbz}) = \sum_{\nu, \mu}c_{\nu, \mu} \bbz^\nu \bar{\bbz}^\mu$ is called 
{\em polar weighted homogeneous} if there exists a weight vector 
$Q={}^t(q_1, \dots , q_n)\ (\neq \zero) \in \bz^n$ 
and 
an integer $d_p$ ($>0,\ 0\ \text{or}\ <0$) such that 
$$
c_{\nu, \mu} \neq 0 \Longrightarrow Q(\nu - \mu) = \sum_{i=1}^n q_i(\nu_i - \mu_i) = d_p .
$$
We call $d_p$ {\em the polar degree} of $f$, and define 
$
\pdeg_Q f := d_p .
$
\end{itemize}
}
\end{definition}

Every face function $f_{\Delta(P)}(\bbz)$, where $P$ is strictly positive, 
of a mixed function germ $(f,\zero)$ 
is a radially weighted homogeneous polynomial of positive radial degree $d(P)$ with respect to the weight vector $P$. 

\begin{example}[\cite{Oka2018}, Example 9.17]
$f(\bbz, \bar{\bbz}) := z_1^2 \bar{z_1} - z_2 \bar{z_2}^2$ is radially weighted homogeneous with respect to $P={}^t(1,1)$ and 
polar weighted homogeneous with respect to $Q={}^t(1,-1)$. 
\end{example}

\medskip

For radially and polar weighted homogeneous polynomials, we have the following basic facts: 
\begin{lemma}[cf.\cite{Oka2008}, \cite{Oka2010}]\label{mixed-action}
Let $f(\bbz, \bar{\bbz})$ be a mixed polynomial. We have the following. 
\begin{itemize}
\item Let $P={}^t(p_1, \dots , p_n)\ (\neq \zero) \in \bz_+^n$ be a weight vector and $d_r$ be a positive integer. 
For a positive real number $t$ and $\bbz \in \bc^n$, we define 
$$
t \circ \bbz := (t^{p_1} z_1, \dots , t^{p_n} z_n) .   
$$
Then, 
$
f(t \circ \bbz) = t^{d_r}f(\bbz)
$
for every $t>0$ and every $\bbz \in \bc^n$ 
if and only if 
``$c_{\nu, \mu} \neq 0 \Rightarrow P(\nu + \mu) = d_r$" holds. 

\item Let $Q={}^t(q_1, \dots , q_n)\ (\neq \zero) \in \br^n$ be a weight vector and $d_p$ be an integer. 
For a real number $\theta$ and $\bbz \in \bc^n$, we define 
$$
\theta \circ \bbz := (e^{iq_1 \theta} z_1, \dots , e^{iq_n \theta} z_n) .  
$$
Then, 
$
f(\theta \circ \bbz) = e^{id_p \theta}f(\bbz)
$
for every $\theta \in \br$ and every $\bbz \in \bc^n$ 
if and only if 
``$c_{\nu, \mu} \neq 0 \Rightarrow Q(\nu - \mu) = d_p$" holds. 
\hfill \QED
\end{itemize}
\end{lemma}

\begin{proposition}[Euler equalities,\ \cite{Oka2018}]\label{Euler}
Let $f(\bbz,\overline{\bbz})=\sum_{\nu,\mu}c_{\nu,\mu}\bbz^{\nu}\overline{\bbz}^{\mu}$ be a mixed polynomial. 
\begin{description}
\item[(R)]\ If $f(\bbz,\overline{\bbz})$ is a radially weighted homogeneous polynomial 
of radial degree $d_{r}\ (>0)$ with respect to a weight vector $P={}^{t}\!(p_{1},\dots,p_{n})$, then we have 
\begin{equation}\label{reulereq}
\sum_{j=1}^{n}p_{j}\left(z_{j}\frac{\partial f}{\partial z_{j}} + \overline{z_{j}}\frac{\partial f}{\partial \overline{z_{j}}} \right)
= 
d_{r}f(\bbz,\overline{\bbz}) .
\end{equation}
\item[(P)]\ If $f(\bbz,\overline{\bbz})$ is a polar weighted homogeneous polynomial 
of polar degree $d_{p}$ with respect to a weight vector $Q={}^{t}\!(q_{1},\dots,q_{n})$, then we have 
\begin{equation}\label{peulereq}
\sum_{j=1}^{n}q_{j}\left(z_{j}\frac{\partial f}{\partial z_{j}} - \overline{z_{j}}\frac{\partial f}{\partial \overline{z_{j}}} \right)
= 
d_{p}f(\bbz,\overline{\bbz}) . \ \ \Box
\end{equation}
\end{description}
\end{proposition}

\medskip

\begin{definition}[mixed weighted homogeneous polynomial,\ \cite{Oka2018}, pp.182--184]\label{mixed-whp}  
{\rm 
We say a mixed polynomial $f(\bbz, \bar{\bbz})$ is 
a {\em mixed weighted homogeneous polynomial}
if 
it is both radially and polar weighted homogeneous. 
Here, the corresponding weight vectors $P$ and $Q$ are possibly different. 
We say a mixed weighted homogeneous polynomial $f(\bbz, \bar{\bbz})$ is 
a {\em strongly mixed weighted homogeneous polynomial} if 
$f$ is radially and polar weighted homogeneous 
with respect to the \underline{same} weight vector $P$. 
Furthermore, 
a mixed weighted homogeneous polynomial $f$ is called 
a {\em strongly polar positive mixed weighted homogeneous polynomial} 
with respect to a weight vector $P$ 
if 
$f$ is radially and polar weighted homogeneous with respect to the same weight vector $P$ and 
$\pdeg_P f >0$. 
}
\end{definition}

\begin{definition}[cf.\cite{Oka2018},Definition 9.18;\ \cite{Oka2015},p.174]\label{strongly-mixed-wh-face-type}  
{\rm 
Let $(f,\zero)$ be a mixed function germ at $\zero \in \bc^n$. 
The germ $(f,\zero)$ of a mixed function $f(\bbz, \bar{\bbz})$ at $\zero \in \bc^n$ is defined to be 
{\em of strongly polar positive mixed weighted homogeneous face type} if 
for every compact face 
\footnote{
In \cite{Oka2015}, a convenient (Definition \ref{convenient-function}) mixed function germ $(f(\bbz, \bar{\bbz}),\zero)$ is defined to be 
{\em of strongly polar positive mixed weighted homogeneous face type} if 
the face function $f_\Delta(\bbz, \bar{\bbz})$ is a strongly polar positive mixed weighted homogeneous polynomial 
for every \underline{$(n-1)$-dimensional face}. 
Besed on this definition, 
Proposition 10 of \cite{Oka2015} asserts that 
the face function $f_{\Delta(P)}$ is also a strongly polar positive mixed weighted homogeneous polynomial with respect to $P$ 
for \underline{every weight vector $P$} 
when $(f,\zero)$ is a convenient mixed function germ of strongly polar positive mixed weighted homogeneous face type in the sense of \cite{Oka2015}. 
}
$\Delta$, 
the face function $f_\Delta(\bbz, \bar{\bbz})$ is 
a strongly polar positive mixed weighted homogeneous polynomial (Definition \ref{mixed-whp}) 
with respect to 
\underline{some} strictly positive weight vector $P$ with $\Delta = \Delta(P)$. 
}
\end{definition}

\medskip

\subsection{Newton non-degeneracy and Strong Newton non-degeneracy}\label{section-Newton-non-degeneracy}

\begin{definition}[\cite{Oka2010}, p.6, Definition 3;\ \cite{Oka2018}, p.80 and pp.181--182]\label{Newton-non-degeneracy}
{\rm 
Let $(f,\zero)$ be the germ of a mixed function $f$ at $\zero \in \bc^n$. 
\begin{enumerate}
\item We say $(f,\zero)$ is {\em Newton non-degenerate 
over a compact face $\Delta$} 
if 
$0$ is not a mixed critical value of the face function 
$f_\Delta : {\bc^*}^n \to \bc$. 
(In particular, 
if $f_\Delta^{-1}(0) \cap {\bc^*}^n = \emptyset$, then $0$ is not a mixed critical value of the face function 
$f_\Delta : {\bc^*}^n \to \bc$.) 

\item Let $\Delta$ be a compact face with $\dim \Delta \geq 1$. 
We say $(f,\zero)$ is {\em strongly Newton non-degenerate 
over $\Delta$} 
if 
the face function $f_\Delta : {\bc^*}^n \to \bc$ has no mixed critical points 
and $f_\Delta : {\bc^*}^n \to \bc$ is surjective onto $\bc$. 

\item Let $\Delta$ be a compact face with $\dim \Delta = 0$, that is, 
$\Delta$ is a vertex of $\Gamma_+(f)$. 
We say $(f,\zero)$ is {\em strongly Newton non-degenerate 
over $\Delta$} 
if 
the face function $f_\Delta : {\bc^*}^n \to \bc$ has no mixed critical points. 
\end{enumerate}
}
\end{definition}

\begin{definition}[\cite{Oka2018}, p.80 and p.182]\label{Newton-non-degenerate-germ}
{\rm 
We say the germ $(f,\zero)$ of a mixed function $f$ at $\zero \in \bc^n$ is 
{\em Newton non-degenerate} (respectively, {\em strongly Newton non-degenerate}) 
if 
$(f,\zero)$ is Newton non-degenerate (respectively, strongly Newton non-degenerate) 
over every compact face $\Delta$. 
}
\end{definition}

\begin{example}[\cite{Oka2018},\ 8.3.1]
The germ of the mixed homogeneous polynomial 
$\rho(\bbz,\overline{\bbz}) = \sum_{j=1}^n z_j \overline{z_j} = \sum_{j=1}^{n}|z_{j}|^{2}$  
of degree $2$ at $\zero \in \bc^{n}$
is Newton non-degenerate, but not strongly Newton non-degenerate. 
Actually, for every compact face $\Delta$, we have ${\rho_{\Delta}}^{-1}(0) \cap \bc^{*n} = \emptyset$. 
However, since $\rho_{\Delta}(\bc^{n}) \subset \br$, 
every point in $\bc^{*n}$ is a mixed critical point. Hence, $\rho$ is not strongly Newton non-degenerate. 
\end{example}

Note that if $f^{-1}(0) \cap {\bc^*}^n \neq \emptyset$, then $f : {\bc^*}^n \to \bc$ is surjective. 
Namely, in this case, Newton non-degeneracy over a compact face $\Delta(P)$ implies 
strong Newton non-degeneracy over $\Delta(P)$. 
For more precise statements, 
see Proposition \ref{remark4} below: 
\begin{proposition}[Remark 4 in \cite{Oka2010};\ and also \cite{Saito-Takashimizu2021winter}, \cite{Saito-Takashimizu2022}]\label{remark4}
\ \ 

\begin{itemize}
\item Let $f(\bbz)$ be a holomorphic weighted homogeneous polynomial 
of positive degree with respect to a strictly positive weight vector $P$. 
(Then, $f = f_{\Delta(P)}$.)

{\bf (i)}\ Suppose that $f$ is Newton non-degenerate over $\Delta(P)$, 
namely, 
$0$ is not a critical value of $f_{\Delta(P)}=f : {\bc^*}^n \to \bc$. 
Then $f_{\Delta(P)}=f : {\bc^*}^n \to \bc$ has no critical point. 
Hence, with {\bf (ii)} below, $f$ is strongly Newton non-degenerate over $\Delta(P)$. 

{\bf (ii)}\ Suppose that $\dim \Delta(P) \geq 1$, namely, $f = f_{\Delta(P)}$ has at least two monomials. 
Then $f_{\Delta(P)}=f : {\bc^*}^n \to \bc$ is surjective. 

\item Let $f(\bbz, \bar{\bbz})$ be a mixed weighted homogeneous polynomial (Definition \ref{mixed-whp}) with respect to 
a radial weight vector $P\ (\gg 0)$ and a polar weight vector $Q$. 

{\bf (iii)}\ Suppose that $(f,\zero)$ is Newton non-degenerate 
\underline{over a compact face $\Delta(P)$} 
and 
the face function $f_{\Delta(P)}=f$ is a polar weighted homogeneous polynomial of \underline{non-zero} polar degree  
with respect to the polar weight vector $Q$. 
Then $f : {\bc^*}^n \to \bc$ has no mixed critical points. 

{\bf (iv)}\ In addition to {\bf (iii)}, 
we assume that 
$f^{-1}(0) \cap {\bc^*}^n \neq \emptyset$.  
Then $f : {\bc^*}^n \to \bc$ is surjective. 
Hence, with {\bf (iii)}, 
in this case, Newton non-degeneracy \underline{over a compact face $\Delta(P)$} implies 
strong Newton non-degeneracy \underline{over $\Delta(P)$}. \hfill \QED
\end{itemize}
\end{proposition}


\bigskip

\section{Mixed weighted homogeneous polynomials of type $J_{10}^-$} \label{section-J10}

We now consider 
the weighted homogeneous polynomials of type $J_{10}^-$:
\begin{equation}\label{J10-}
\begin{array}{ccl}
f(z_1,z_2) & := & (z_2 - z_1^2)(z_2 - 2z_1^2)(z_2 - k z_1^2)\\
             & =  & z_2^3 - (k +3)z_1^2z_2^2 + (3k +2)z_1^4z_2 - 2k z_1^6 ,
\end{array}
\end{equation}
where we assume that $k \in \br$ and $k>2$. 
Note that these polynomials are weighted homogeneous with respect to the weight vector 
$$P := {}^t(1,2)$$
of degree $6$ with real coefficients. 

Based on the weighted homogeneous polynomials \eqref{J10-}, 
let us consider the radially weighted homogeneous polynomials 
with respect to the radial weight vector $P = {}^t(1,2)$ of radial degree $6$:
\begin{equation}\label{radialJ10-}
f_{a,b,c,d,e,\mathrm{f}} := 
z_2^a \bar{z}_2^{3-a}  - (k +3)z_1^b \bar{z}_1^{2-b} z_2^c \bar{z}_2^{2-c} + (3k +2)z_1^d \bar{z}_1^{4-d} z_2^e \bar{z}_2^{1-e}  - 2k z_1^\mathrm{f} \bar{z}_1^{6-\mathrm{f}},
\end{equation}
where $a, b, c, d, e, \mathrm{f}$ are integers with $0\leq a \leq 3,\ 0\leq b,c \leq 2,\ 0\leq d \leq 4,\ 0\leq e \leq 1$ and $0\leq \mathrm{f} \leq 6$. 

\medskip

Obviously $f_{a,b,c,d,e,\mathrm{f}}$ are convenient (Definition \ref{convenient-function}). 

\subsection{The radial Newton polyhedron, the dual Newton diagram and the regular simplicial cone subdivision}

The radial Newton polyhedron of the mixed function germ $f_{a,b,c,d,e,\mathrm{f}}$ at $\zero$ is as follows:
\begin{figure}[H]
\begin{center}
\begin{picture}(120,120)
\multiput(16,0)(16,0){7}{\line(0,1){78}}
\multiput(0,16)(0,16){4}{\line(1,0){124}}
\put(0,0){\line(0,1){47}}
\put(0,0){\line(1,0){93}}
\put( -2,-10){{\tiny $0$}}
\put( 14,-10){{\tiny $1$}}
\put( 30,-10){{\tiny $2$}}
\put( 46,-10){{\tiny $3$}}
\put( 62,-10){{\tiny $4$}}
\put( 77,-10){{\tiny $5$}}
\put( 93,-10){{\tiny $6$}}
\put(109,-10){{\tiny $7$}}
\put(-8, -1){{\tiny $0$}}
\put(-8, 15){{\tiny $1$}}
\put(-8, 31){{\tiny $2$}}
\put(-8, 47){{\tiny $3$}}
\put(-8, 63){{\tiny $4$}}
\put(-14,84){{\footnotesize $\nu_2 + \mu_2$}}
\put(131, -2){{\footnotesize $\nu_1 + \mu_1$}} 

\put(96, 0){\circle*{4}}  
\put(64,16){\circle*{4}}  
\put(32,32){\circle*{4}}  
\put( 0,48){\circle*{4}}  

\linethickness{0.3mm}    
\put(0,48){\vector(0,1){32}}
\put(96,0){\vector(1,0){30}}
\multiput(0,48)(2,-1){4}{\line(2,-1){91}}   
\end{picture}

\vspace{8mm}

\end{center}
\caption{The radial Newton polyhedron of $f_{a,b,c,d,e,\mathrm{f}}$}
\end{figure}
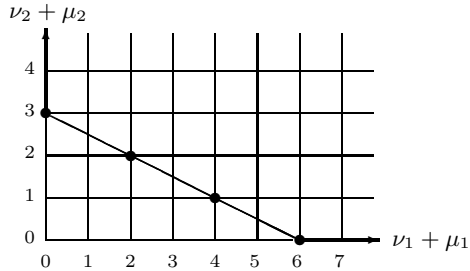

\bigskip

\begin{figure}[H]
\begin{center}
\includegraphics[width=4.5cm]{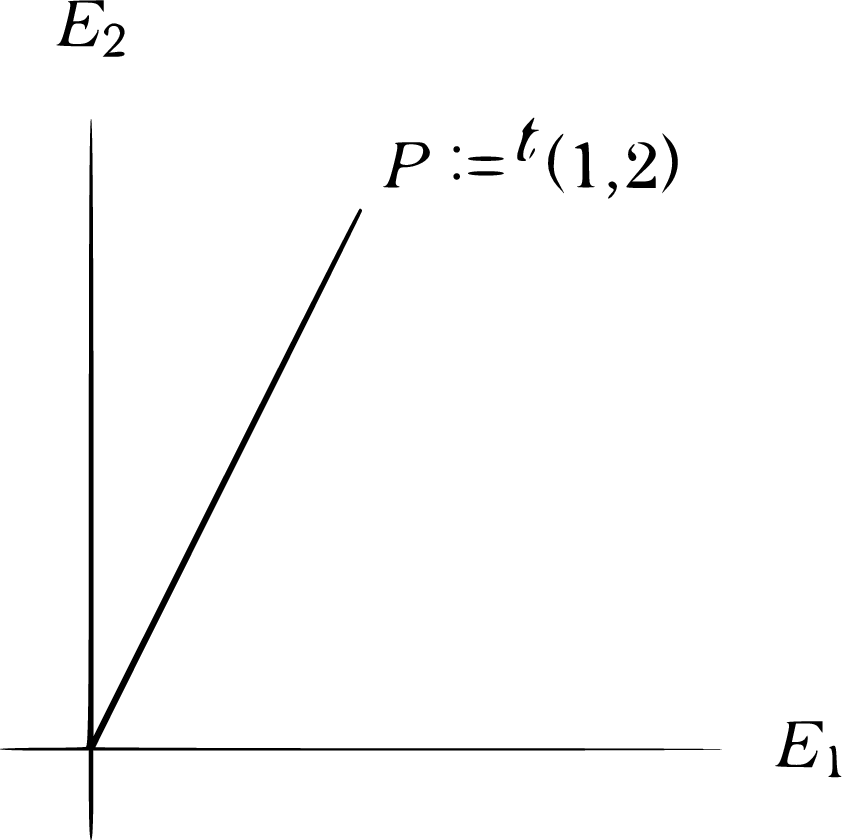}
\end{center}
\caption{The dual Newton diagram $\Gamma^*(f)$}
\label{dual_newton_diagram}
\end{figure}

The dual Newton diagram $\Gamma^*(f)$ of the mixed function germ $f := f_{a,b,c,d,e,\mathrm{f}}$ at $\zero$ is 
as in the Figure $\ref{dual_newton_diagram}$, 
where we set 
$E_1 = {}^t(1,0),\ P= {}^t(1,2),\ \ E_2 = {}^t(0,1)$. 

Adding the vertex $S := \ ^t(1,1)$ to $\Gamma^*(f)$, 
we obtain a regular simplicial cone subdivision (see \cite{Oka2018}, \S 5.6; or \cite{Oka1997} for the definition) 
$\Sigma^*$ of $N^+_\br$ 
which is admissible for $\Gamma^*(f)$. 

\begin{figure}[H]
\begin{center}
\includegraphics[width=4.5cm]{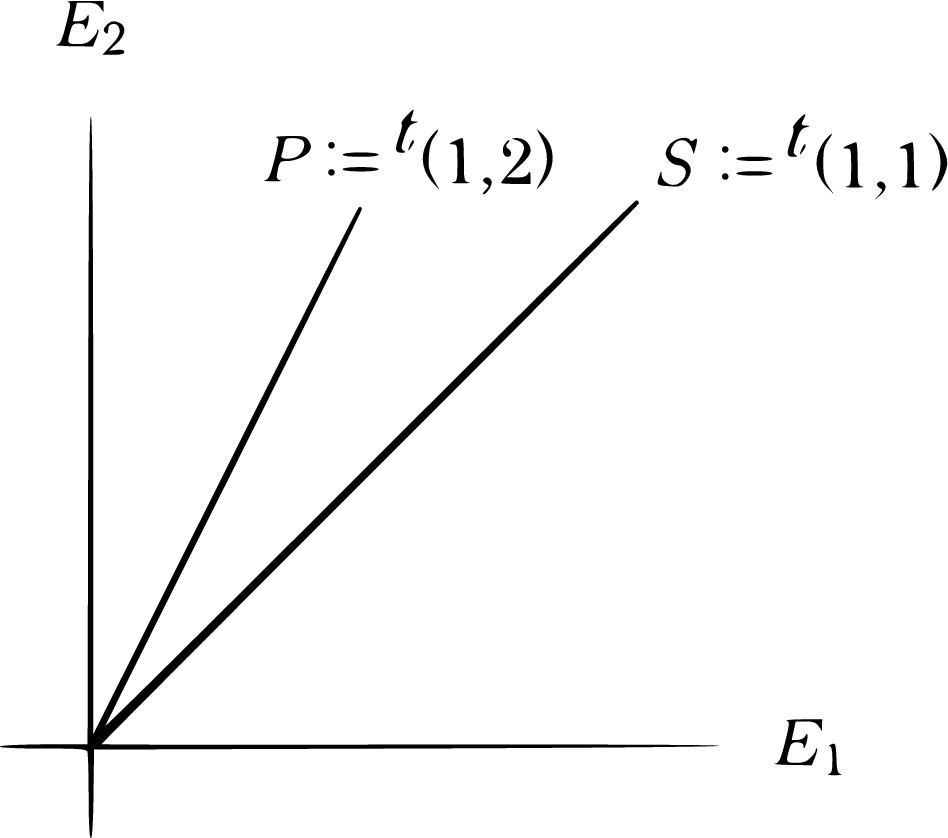}
\end{center}

\caption{The regular simplicial subdivision $\Sigma^*$ which is admissible for $\Gamma^*(f)$.}
\label{regular_subdivision_abcdef}
\end{figure}

\medskip

\begin{definition}[\cite{Oka2018}, p.98] \label{convenient-subdivision-f}  
{\rm 
We set 
$E_j := {}^t(\underbrace{0, \dots , 0, 1}_{j}, 0, \dots , 0) \ (\in \bz^n)$ 
for every $j \ (= 1, \dots , n)$, and 
$$E_J := \Cone (E_{j_1}, \dots , E_{j_k}),$$
where $J := \{ j_1, \dots , j_k \} \ (\subsetneqq \{1, \dots, n\})$. 
Let $(f,\zero)$ be the germ of a mixed function $f$ at $\zero \in \bc^n$. 
For $I \subset \{1, \dots, n\}$, 
$f^I$ denotes the restriction of $f$ on the coordinate subspace 
$\bc^I := \{ \bbz \,|\, z_j=0,\, j\notin I \}$. 
Let $\Sigma^*$ be a regular simplicial cone subdivision 
which is admissible for $\Gamma^*(f)$. 
We say $\Sigma^*$ is {\em convenient} 
if for every $I$ with $f^I \not \equiv 0$, the cone $E_{I^c}$ is contained in $\Sigma^*$. 
}
\end{definition}

Then our $\Sigma^*$ (Figure \ref{regular_subdivision_abcdef}) is convenient in the sense of Definition \ref{convenient-subdivision-f}, 
namely, 
the cone $E_{I^c}$ is contained in $\Sigma^*$ for every $I$ with $f^I \not \equiv 0$. 

\bigskip

We moreover set $T := {}^t(1,3)$. 
Then, all the \underline{compact} faces of the radial Newton polyhedron of the germ $(f_{a,b,c,d,e,\mathrm{f}},\zero)$ are 
$$\Delta(P),\ \ \Delta(S),\ \text{and} \ \Delta(T),$$
and hence, the face functions of $(f := f_{a,b,c,d,e,\mathrm{f}},\zero)$ are 
$$
f_P = f_{a,b,c,d,e,\mathrm{f}},\ f_S = z_2^a \bar{z}_2^{3-a}\ \text{and} \ f_T = - 2k z_1^\mathrm{f} \bar{z}_1^{6-\mathrm{f}} .
$$

\medskip

\subsection{Strongly mixed weighted homogeneous polynomials of type $J_{10}^-$}

The radially weighted homogeneous polynomials $f_{a,b,c,d,e,\mathrm{f}}$ are 
\underline{strongly} mixed weighted homogeneous with respect to $P$ (Definition \ref{mixed-whp}) 
if and only if 
$(a, b, c, d, e, \mathrm{f})$ are in the following $5$ cases:

\begin{center}
\begin{tabular}{c||cccccc|c|l}\label{strongly-mixed-whp-5cases}
    & $a$ & $b$ & $c$ & $d$ & $e$ & $\mathrm{f}$ & {\small polar degree}  &  \\
\hline
I   & 3 & 2 & 2 & 4 & 1 & 6 &  6   & holomorphic case \\
II  & 2 & 2 & 1 & 4 & 0 & 4 &  2   &  \\
III & 2 & 0 & 2 & 4 & 0 & 4 &  2   &  \\
IV & 2 & 2 & 1 & 2 & 1 & 4 &  2   &  \\
V  & 2 & 0 & 2 & 2 & 1 & 4 &  2  &  
\end{tabular}
\end{center}

In the cases II $\sim$ V ($a=2,\ \mathrm{f}=4$), the radial degrees of $f := f_{2,b,c,d,e,4}$ are $6$, and 
the polar degrees of those are $2$. 
For their $0$-dimensional face functions 
$f_S = z_2^2 \bar{z}_2$ and $f_T = - 2k z_1^4 \bar{z}_1^2$, we have 
$$
\rdeg_S f_S = 1\cdot 0 + 1\cdot (2+1) = 3
, \ \ \ 
\pdeg_S f_S = 1\cdot 0 + 1\cdot (2-1) = 1, 
$$
$$
\rdeg_T f_T = 1\cdot (4+2) + 3\cdot 0 = 6
, \ \text{and} \ 
\pdeg_T f_T = 1\cdot (4-2) + 3\cdot 0 = 2 . 
$$
Hence, remark that the germs $(f_{2,b,c,d,e,4},\zero)$ (in the cases II $\sim$ V) 
are of 
strongly polar positive mixed weighted homogeneous face type (Definition \ref{strongly-mixed-wh-face-type}). 

\medskip

\subsection{Strong Newton non-degeneracy}

We now argue about the Newton degeneracy of 
the strongly mixed weighted homogeneous polynomials $f_{2,b,c,d,e,4}$ in the cases II $\sim$ V. 
Especially let us investigate the case IV with $k = 3$:
\begin{equation}\label{IV-alpha3}
f_{2,2,1,2,1,4} \ (k=3) = 
z_2^2 \bar{z}_2  - 6z_1^2 z_2 \bar{z}_2 + 11z_1^2 \bar{z}_1^2 z_2 - 6 z_1^4 \bar{z}_1^2 .
\end{equation}

\begin{lemma}\label{no-mixed-critical-points}
We have the following. 
\begin{enumerate}
\item $f := f_{2,2,1,2,1,4} :  \ (k=3)$ has no mixed critical points on ${\bc^*}^2$. 
\item $f : {\bc^*}^2 \to \bc$ is surjective. 
\item The $0$-dimensional face functions $f_S,\ f_T$ of $(f,\zero)$ also have has no mixed critical points on ${\bc^*}^2$. 
\end{enumerate}
Thus, the germ $(f,\zero)$ is strongly Newton non-degenerate (Definition \ref{Newton-non-degenerate-germ}). 
\end{lemma}

\begin{proof}
We first prove the assertion (1). 
We have 
\begin{equation}\label{partial}
\begin{array}{ccl}
\vspace{0.3cm}
\displaystyle \frac{\partial f}{\partial z_1}         & = & -12z_1 z_2 \bar{z}_2 + 22z_1 \bar{z}_1^2 z_2 - 24 z_1^3 \bar{z}_1^2 ,\\
\displaystyle \frac{\partial f}{\partial \bar{z}_1} & = & 22 z_1^2 \bar{z}_1 z_2 - 12 z_1^4 \bar{z}_1 \\
                                               & = & 2z_1^2 \bar{z}_1 (11z_2 - 6z_1^2) ,\\
\displaystyle \frac{\partial f}{\partial z_2}         & = & 2z_2 \bar{z}_2  - 6z_1^2 \bar{z}_2  + 11z_1^2 \bar{z}_1^2\\
                                               & = & 2|z_2|^2 + z_1^2(11 \bar{z}_1^2 - 6 \bar{z}_2) ,\\
\displaystyle \frac{\partial f}{\partial \bar{z}_2} & = & z_2^2   - 6z_1^2 z_2 \\
                                               & = & z_2(z_2 - 6z_1^2) .
\end{array}
\end{equation}
By Lemma \ref{mixed-action}, we have 
$
f(tz_1, t^2z_2) = r^6 e^{2i\theta} f(z_1, z_2) 
$
for every $t = re^{i\theta}\ (\in \bc^*), r>0, \theta \in \br$ 
since $d_r = 6$ (radial degree) and $d_p = 2$ (polar degree). 
Hence, if $f(z_1, z_2)=0$, then we have 
\begin{equation}\label{t-action}
f(tz_1, t^2z_2)=0
\end{equation}
for all $t \in \bc^*$. 
By Proposition \ref{Euler}, we have 
$$
\frac{6+2}{2}f(z_1, z_2) = z_1 \frac{\partial f}{\partial z_1} + 2z_2 \frac{\partial f}{\partial z_2}
$$
and 
$$
\frac{6-2}{2}f(z_1, z_2) = \bar{z}_1 \frac{\partial f}{\partial \bar{z}_1} + 2 \bar{z}_2 \frac{\partial f}{\partial \bar{z}_2} .
$$

Suppose that $\ba = (a_1, a_2) \in {\bc^*}^2$ is a mixed critical point of $f$. 
Recall that the following two conditions are equivalent (Proposition 1 in \cite{Oka2008}):
\begin{enumerate}
\item $\ba = (a_1, a_2)\ (\in \bc^2)$ is a mixed critical point of $f$. 
\item There exists a complex number $\alpha$ with $|\alpha|=1$ which satisfies 
$$\left(\overline{\frac{\partial f}{z_1}(\ba)},\ \overline{\frac{\partial f}{z_2}(\ba)} \right) 
= \alpha \left(\frac{\partial f}{\partial \bar{z_1}}(\ba),\ \frac{\partial f}{\partial \bar{z_2}}(\ba) \right) .$$
\end{enumerate}

\medskip

\noindent
We have 
$$
\begin{array}{ccl}
\vspace{0.2cm}
4\overline{f(a_1, a_2)} & = & \bar{a}_1 \overline{\frac{\partial f}{\partial z_1}(a_1, a_2)} + 2\bar{a}_2 \overline{\frac{\partial f}{\partial z_2}(a_1, a_2)}\\
\vspace{0.2cm}
                              & = & \bar{a}_1 \alpha \frac{\partial f}{\partial \bar{z_1}}(a_1, a_2) + 2\bar{a}_2 \alpha \frac{\partial f}{\partial \bar{z_2}}(a_1, a_2)\\
                              & = & \alpha (\bar{a}_1 \frac{\partial f}{\partial \bar{z_1}}(a_1, a_2) + 2\bar{a}_2 \frac{\partial f}{\partial \bar{z_2}}(a_1, a_2))
\end{array}
$$
and 
$$
2f(a_1, a_2) = \bar{a}_1 \frac{\partial f}{\partial \bar{z}_1}(a_1, a_2) + 2 \bar{a}_2 \frac{\partial f}{\partial \bar{z}_2}(a_1, a_2) .
$$
Thus we have 
$
4\overline{f(a_1, a_2)} = 2\alpha f(a_1, a_2)
$
and 
$
4|f(a_1, a_2)| = 2|\alpha | |f(a_1, a_2)| = 2|f(a_1, a_2)| .
$
Hence, it is concluded that 
$f(a_1, a_2)=0$. 
Thus we see that 
every mixed critical point of $f$ is a zero of that. 

\medskip

Now we show that $f$ has no mixed critical points on ${\bc^*}^2$. 
It is sufficient to prove that 
$f$ has no mixed critical points in $f^{-1}(0) \cap {\bc^*}^2$, 
namely, 
$f$ is Newton non-degenerate \underline{over $\Delta(P)$}. 

Suppose that $(a_1, a_2) \in  {\bc^*}^2$ and $f(a_1, a_2)=0$. 
If 
$|\frac{\partial f}{\partial z_1}(a_1, a_2)| \neq |\frac{\partial f}{\partial \bar{z_1}}(a_1, a_2)|$ or 
$|\frac{\partial f}{\partial z_2}(a_1, a_2)| \neq |\frac{\partial f}{\partial \bar{z_2}}(a_1, a_2)|$, 
then $(a_1, a_2)$ is not a mixed critical point. 
Thus we may moreover suppose that 
$|\frac{\partial f}{\partial z_1}(a_1, a_2)| = |\frac{\partial f}{\partial \bar{z_1}}(a_1, a_2)|$ and 
$|\frac{\partial f}{\partial z_2}(a_1, a_2)| = |\frac{\partial f}{\partial \bar{z_2}}(a_1, a_2)|$. 
Then we have 
\begin{equation}\label{alpha_1alpha_2}
\overline{\frac{\partial f}{\partial z_1}(a_1, a_2)}  =  \alpha_1 \frac{\partial f}{\partial \bar{z_1}}(a_1, a_2)\ \ \ \text{and} \ \ \ 
\overline{\frac{\partial f}{\partial z_2}(a_1, a_2)}  =  \alpha_2 \frac{\partial f}{\partial \bar{z_2}}(a_1, a_2)
\end{equation}
for some $\alpha_1, \alpha_2 \in \bc$ with $|\alpha_1| = |\alpha_2| = 1$. 
On the other hand, by \eqref{t-action}, we have 
$
f(ta_1, t^2a_2)=0
$
for all $t \in \bc^*$. 
We set $\mathbf{c}(t)=(c_1(t),c_2(t)):=(ta_1,t^2a_2)$. 
Then, by the chain rules for Wirtinger derivative, we have 
$$
\begin{array}{ccl}
\vspace{0.2cm}
\frac{d}{dt}(f \circ \mathbf{c})(t) 
& = & \sum_{j=1}^2(\frac{\partial f}{\partial z_j} \circ \mathbf{c})(t)\frac{dc_j}{dt}(t) + \sum_{j=1}^2(\frac{\partial f}{\partial \bar{z}_j} \circ \mathbf{c})(t)\frac{d\overline{c_j}}{dt}(t) \\
& = & \sum_{j=1}^2(\frac{\partial f}{\partial z_j} \circ \mathbf{c})(t)\frac{dc_j}{dt}(t) .
\end{array}
$$
Hence, we have 
$$
\frac{\partial f}{\partial z_1}(ta_1,t^2a_2)\frac{dc_1}{dt}(t) + \frac{\partial f}{\partial z_2}(ta_1,t^2a_2)\frac{dc_2}{dt}(t) =0.
$$
Setting $t=1$, we have
$$
a_1 \frac{\partial f}{\partial z_1}(a_1, a_2) + 2 a_2 \frac{\partial f}{\partial z_2}(a_1, a_2) = 0. 
$$
By \eqref{alpha_1alpha_2}, we have 
\begin{equation}\label{chain-rule1}
\bar{a}_1 \alpha_1 \frac{\partial f}{\partial \bar{z_1}}(a_1, a_2) + 2 \bar{a}_2 \alpha_2 \frac{\partial f}{\partial \bar{z_2}}(a_1, a_2) = 0. 
\end{equation}
We now show that such $(a_1, a_2)$ is \underline{not} a mixed critical point. 
It is sufficient to prove that $\alpha_1 \neq \alpha_2$. 
Suppose that $\alpha_1 = \alpha_2$. 
We have 
\begin{equation}\label{chain-rule2}
\bar{a}_1 \frac{\partial f}{\partial \bar{z_1}}(a_1, a_2) + 2 \bar{a}_2 \frac{\partial f}{\partial \bar{z_2}}(a_1, a_2) = 0, 
\end{equation}
that is, 
$$
\bar{a}_1 a_1^2 \bar{a}_1 (11a_2 - 6a_1^2) + \bar{a}_2 a_2(a_2 - 6a_1^2) = 0 .
$$
We have 
\begin{equation}\label{11-6_1-6}
|a_1|^4 (11a_2 - 6a_1^2) + |a_2|^2(a_2 - 6a_1^2) = 0 ,
\end{equation}
$$
11|a_1|^4 a_2 - 6 |a_1|^4 a_1^2 + |a_2|^2 a_2 - 6 |a_2|^2 a_1^2 = 0 , 
$$
$$
(11|a_1|^4  + |a_2|^2) a_2  =  6 (|a_1|^4  + |a_2|^2) a_1^2 , 
$$
and
\begin{equation}\label{a_2a_1sq2}
a_2  =  \frac{6(|a_1|^4  + |a_2|^2)}{11|a_1|^4  + |a_2|^2}\, a_1^2 .
\end{equation}
This means that $a_2$ and $a_1^2$ have the same direction in the complex plane $\bc$. 
%
%
Note that by \eqref{11-6_1-6}, we also have 
\begin{equation}\label{11-6}
11a_2 - 6a_1^2 = - \frac{|a_2|^2}{|a_1|^4}(a_2 - 6a_1^2) 
\end{equation}
and
\begin{equation}\label{1-6}
a_2 - 6a_1^2 = \frac{-10|a_1|^4}{|a_1|^4 + |a_2|^2}\, a_2 .
\end{equation}
By the calculations \eqref{partial} of partial derivatives, we have
$$
\frac{\partial f}{\partial \bar{z}_2}(a_1, a_2) = a_2(a_2 - 6a_1^2) =  \frac{-10|a_1|^4}{|a_1|^4 + |a_2|^2} a_2^2 . 
$$
On the other hand, we have 
$$
\begin{array}{ccl}
\frac{\partial f}{\partial z_2}(a_1, a_2) & = & 2a_2 \bar{a}_2  - 6a_1^2 \bar{a}_2  + 11a_1^2 \bar{a}_1^2 \\
                                                   & = & a_2 \bar{a}_2 + a_2 \bar{a}_2  - 6a_1^2 \bar{a}_2  + 11a_1^2 \bar{a}_1^2 \\
\vspace{0.1cm}
                                                   & = & a_2 \bar{a}_2 + (a_2 - 6a_1^2)\bar{a}_2 + 11a_1^2 \bar{a}_1^2 \\
\vspace{0.1cm}
                                                   & = & a_2 \bar{a}_2 + \frac{-10|a_1|^4}{|a_1|^4 + |a_2|^2} a_2\bar{a}_2 + 11a_1^2 \bar{a}_1^2 \\
\vspace{0.1cm}
                                                  & = & |a_2|^2 + \frac{-10|a_1|^4}{|a_1|^4 + |a_2|^2} |a_2|^2 + 11|a_1|^4 \\
\vspace{0.1cm}
                                                  & = & \frac{-10|a_1|^4 + |a_1|^4 + |a_2|^2}{|a_1|^4 + |a_2|^2} |a_2|^2 + 11|a_1|^4 \\
\vspace{0.1cm}
                                                 & = & \frac{-9|a_1|^4 + |a_2|^2}{|a_1|^4 + |a_2|^2} |a_2|^2 + 11|a_1|^4 \\
\vspace{0.1cm}
                                                  & = & \frac{-9|a_1|^4 |a_2|^2 + |a_2|^4 + 11|a_1|^8 + 11|a_1|^4 |a_2|^2}{|a_1|^4 + |a_2|^2} \\
\vspace{0.1cm}
                                                  & = & \frac{2|a_1|^4 |a_2|^2 + |a_2|^4 + 11|a_1|^8}{|a_1|^4 + |a_2|^2} > 0 .
\end{array}
$$
By the latter equality of \eqref{alpha_1alpha_2}, we have 
$$
\frac{2|a_1|^4 |a_2|^2 + |a_2|^4 + 11|a_1|^8}{|a_1|^4 + |a_2|^2} = \frac{-10|a_1|^4}{|a_1|^4 + |a_2|^2}\, \alpha_2 a_2^2 .
$$
Since $\alpha_2 a_2^2$ is a negative real number, we have $\alpha_2 \frac{a_2^2}{|a_2|^2} = -1$. 
Namely, we have 
$$
\alpha_1  = \alpha_2  = -\frac{|a_2|^2}{a_2^2} .
$$
Then, on the other hand, we have 
$$
\overline{\frac{\partial f}{\partial z_1}(a_1, a_2)} = -12\bar{a}_1 \bar{a}_2 a_2 + 22 \bar{a}_1 a_1^2 \bar{a}_2 - 24 \bar{a}_1^3 a_1^2
$$
and
$$
\begin{array}{ccl}
\alpha_1 \frac{\partial f}{\partial \bar{z}_1}(a_1, a_2) & = & 2 \alpha_1 a_1^2 \bar{a}_1 (11a_2 - 6a_1^2) \\
\vspace{0.1cm}
 & = & -2 \frac{|a_2|^2}{a_2^2} a_1^2 \bar{a}_1 (11a_2 - 6a_1^2) \\
 & = & -2 \frac{\bar{a}_2}{a_2} a_1^2 \bar{a}_1 (11a_2 - 6a_1^2) .
\end{array}
$$
Thus, by the first equality of \eqref{alpha_1alpha_2}, we have 
$$
-12\bar{a}_1 \bar{a}_2 a_2 + 22 \bar{a}_1 a_1^2 \bar{a}_2 - 24 \bar{a}_1^3 a_1^2
=  -2 \frac{\bar{a}_2}{a_2} a_1^2 \bar{a}_1 (11a_2 - 6a_1^2) , 
$$
$$
6\bar{a}_1 \bar{a}_2 a_2 - 11 \bar{a}_1 a_1^2 \bar{a}_2 + 12 \bar{a}_1^3 a_1^2
=  \frac{\bar{a}_2}{a_2} a_1^2 \bar{a}_1 (11a_2 - 6a_1^2) , 
$$
and 
\begin{equation}\label{6-11+12}
6 \bar{a}_2 a_2 - 11 a_1^2 \bar{a}_2 + 12 \bar{a}_1^2 a_1^2
=  \frac{\bar{a}_2}{a_2} a_1^2 (11a_2 - 6a_1^2) .
\end{equation}
%
%
%
The equality \eqref{6-11+12} with \eqref{11-6} and \eqref{1-6} yields 
$$
\begin{array}{ccl}
\vspace{0.1cm}
6 \bar{a}_2 a_2 - 11 a_1^2 \bar{a}_2 + 12 \bar{a}_1^2 a_1^2 & = &  - \frac{\bar{a}_2}{a_2} a_1^2 \frac{|a_2|^2}{|a_1|^4}(a_2 - 6a_1^2)\\
\vspace{0.1cm}
                                                                                  & = & \frac{\bar{a}_2}{a_2} a_1^2 \frac{|a_2|^2}{|a_1|^4} \frac{10|a_1|^4}{|a_1|^4 + |a_2|^2} a_2\\
                                                                                  & = & \frac{10|a_2|^2}{|a_1|^4 + |a_2|^2} a_1^2 \bar{a}_2 .
\end{array}
$$
Thus we have 
$$
12 |a_1|^4 + 6 |a_2|^2 = \left( \frac{10|a_2|^2}{|a_1|^4 + |a_2|^2} + 11 \right) a_1^2 \bar{a}_2 .
$$
Using \eqref{a_2a_1sq2}, we have 
\begin{equation}
\begin{array}{ccl}
\vspace{0.1cm}
12 |a_1|^4 + 6 |a_2|^2 & = & \left( \frac{10|a_2|^2 + 11|a_1|^4 + 11|a_2|^2}{|a_1|^4 + |a_2|^2} \right) \frac{6(|a_1|^4  + |a_2|^2)}{11|a_1|^4  + |a_2|^2}\, a_1^2  \bar{a}_1^2 \\
\vspace{0.1cm}
                              & = & \left( \frac{11|a_1|^4 + 21|a_2|^2}{|a_1|^4 + |a_2|^2} \right) \frac{6(|a_1|^4  + |a_2|^2)}{11|a_1|^4  + |a_2|^2}\, |a_1|^4 \\
                              & = & 6\frac{11|a_1|^4 + 21|a_2|^2}{11|a_1|^4  + |a_2|^2}\, |a_1|^4 ,
\end{array}
\end{equation}
namely, 
$$
\begin{array}{ccl}
\vspace{0.1cm}
2 |a_1|^4 + |a_2|^2 & = & \frac{11|a_1|^4 + 21|a_2|^2}{11|a_1|^4  + |a_2|^2}\, |a_1|^4 \\
\vspace{0.1cm}
                         & = & \left( 1 + \frac{20|a_2|^2}{11|a_1|^4  + |a_2|^2} \right) \, |a_1|^4 \\
                         & = & |a_1|^4 + \frac{20|a_1|^4 |a_2|^2}{11|a_1|^4  + |a_2|^2} .
\end{array}
$$
Hence, we have 
$$
|a_1|^4 + |a_2|^2 = \frac{20|a_1|^4 |a_2|^2}{11|a_1|^4  + |a_2|^2}, 
$$
$$
(11|a_1|^4  + |a_2|^2)(|a_1|^4 + |a_2|^2) = 20|a_1|^4 |a_2|^2, 
$$
$$
11|a_1|^8 + 12|a_1|^4 |a_2|^2 + |a_2|^4 = 20|a_1|^4 |a_2|^2, 
$$
$$
|a_2|^4 - 8|a_1|^4 |a_2|^2 + 11|a_1|^8= 0,    
$$
and 
$$
(|a_2|^2 - (4 - \sqrt{5})|a_1|^4)(|a_2|^2 - (4 + \sqrt{5})|a_1|^4) = 0 .
$$
Hence, we have 
$$
|a_2|^2 = (4 \pm \sqrt{5})|a_1|^4, 
$$
that is, 
$$
|a_2| = \sqrt{4 \pm \sqrt{5}}\, |a_1|^2 .
$$
Since $a_2$ and $a_1^2$ have the same direction in $\bc$ (see \eqref{a_2a_1sq2}), we have 
$$
a_2 = \sqrt{4 \pm \sqrt{5}}\, a_1^2 .
$$
Recall that $(a_1, a_2)$ is a zero of $f$. We have 
$$
\begin{array}{ccl}
0 = f(a_1, a_2) & = & a_2^2 \bar{a}_2  - 6a_1^2 a_2 \bar{a}_2 + 11a_1^2 \bar{a}_1^2 a_2 - 6 a_1^4 \bar{a}_1^2 \\
              & = & a_2 |a_2|^2 - 6a_1^2 |a_2|^2 + 11a_2|a_1|^4 - 6 a_1^2 |a_1|^4 \\
              & = & (a_2 - 6a_1^2) |a_2|^2 + (11a_2 - 6 a_1^2) |a_1|^4 \\
              & = & (\xi_{\pm} - 6)a_1^2 |a_2|^2 + (11\xi_{\pm} -6)a_1^2 |a_1|^4 \\
              & = & (\xi_{\pm} - 6) \xi_{\pm}^2 a_1^2 |a_1|^4  + (11\xi_{\pm} -6)a_1^2 |a_1|^4 , 
\end{array}
$$
where we set $\xi_{\pm} := \sqrt{4 \pm \sqrt{5}}$. 
Then we have 
$$
(\xi_{\pm} - 6) \xi_{\pm}^2 + (11\xi_{\pm} -6) = 0, 
$$
$$
\xi_{\pm}^3 - 6\xi_{\pm}^2 + 11\xi_{\pm} -6 = 0, 
$$
and 
$$
(\xi_{\pm} -1)(\xi_{\pm} -2)(\xi_{\pm} -3) = 0 .
$$
Finally we have 
$$
\xi_{\pm} = 1,\ 2 \ \text{or} \ 3 .
$$
This assertion contradicts to $\xi_{\pm} = \sqrt{4 \pm \sqrt{5}}$. 
We finally have $\alpha_1 \neq \alpha_2$, 
and see that 
$(a_1, a_2)$ is not a mixed critical point. 

\medskip

We next prove the assertion (2), i.e., the surjectivity of $f : {\bc^*}^2 \to \bc$. 
For the mixed polynomial \eqref{radialJ10-}, if $z_2$ is a real variable, then we have 
$$
f_{a,b,c,d,e,\mathrm{f}}(1,z_2) = 
z_2^3  - (k +3)z_2^2 + (3k +2)z_2 - 2k .
$$
Since $k > 2$, the equation $f_{a,b,c,d,e,\mathrm{f}}(1,z_2) = 0$ has a real solution $z_2 \ (\neq 0)$. 
Hence we have 
$$f_{a,b,c,d,e,\mathrm{f}}^{-1}(0) \cap {\bc^*}^2 \neq \emptyset.$$
Recall that 
$f_{a,b,c,d,e,\mathrm{f}}$ are strongly mixed weighted homogeneous 
if $(a, b, c, d, e, \mathrm{f})$ is one of the $5$ cases in Table \ref{strongly-mixed-whp-5cases}. 
By the above results and {\bf (iii)},{\bf (iv)} of Proposition \ref{remark4}, 
we conclude that 
%
%
$f_{2,2,1,2,1,4}\ (k=3) : {\bc^*}^2 \to \bc$ is surjective. 
(Hence, 
$f_{2,2,1,2,1,4}\ (k=3)$ is strongly Newton non-degenerate over the $1$-dimensional face $\Delta(P)$, 
where $P = {}^t(1,2)$. )

\bigskip

Finally, we prove the assertion (3). 
Let us consider the $0$-dimensional face functions 
$f_S = z_2^a \bar{z}_2^{3-a}$ and $f_T = - 2k z_1^\mathrm{f} \bar{z}_1^{6-\mathrm{f}}$ of $f_{a,b,c,d,e,\mathrm{f}}$ in the Table \ref{strongly-mixed-whp-5cases}. 
If $a=2$, then we have 
$$
\frac{\partial f_S}{\partial z_2} = 2z_2\bar{z}_2,\ 
\frac{\partial f_S}{\partial \bar{z}_2} = z_2^2 , 
$$
and hence, $f_S : {\bc^*}^2 \to \bc$ has no mixed critical points. 
If $\mathrm{f}=4$, then we have 
$$
\frac{\partial f_T}{\partial z_1} = -8k z_1^3 \bar{z}_1^2,\ 
\frac{\partial f_T}{\partial \bar{z}_1} = - 4k z_1^4 \bar{z}_1 , 
$$
and hence, $f_T : {\bc^*}^2 \to \bc$ has no mixed critical points. 
Hence, the $0$-dimensional face functions $f_S$ and $f_T$ of $f := f_{2,b,c,d,e,4}$ 
are strongly Newton non-degenerate. 
This completes the proof of Lemma \ref{no-mixed-critical-points}. 
\end{proof}

\medskip

\subsection{Toric modifications associated with the regular simplicial cone subdivision $\Sigma^*$}

Now let us consider the toric modification $\hat{\pi} : X \to \bc^2$ 
associated with the regular simplicial cone subdivision $\Sigma^*$ (Figure \ref{regular_subdivision_abcdef}). 
Note that 
$$
\hat{\pi}^{-1}(\zero) = \hat{E}(S) \cup \hat{E}(P).
$$

All $2$-dimensional cones of $\Sigma^*$ are as follows (up to permutations of vertices): 
\begin{equation}\label{3charts}
\begin{array}{l}
\sigma_1 := \Cone (E_1,S) = 
\begin{pmatrix}
1 & 1 \\
0 & 1
\end{pmatrix}
,\ \ 
\sigma_2 := \Cone (S,P) = 
\begin{pmatrix}
1 & 1 \\
1 & 2
\end{pmatrix}
,\ \ 
\sigma_3 := \Cone (P,E_2) = 
\begin{pmatrix}
1 & 0 \\
2 & 1 
\end{pmatrix}.
\end{array}
\end{equation}

\bigskip

\bigskip

We now show that the Assumption (*) in Theorem 32 of \cite{Saito-Takashimizu2021winter} 
is satisfied for the mixed polynomial germ $(f_{2,2,1,2,1,4}\ (k=3),\zero)$ and our regular simplicial cone subdivision $\Sigma^*$. 

\bigskip

\noindent
{\bf (I)}\ \ 
We first set 
$$
\sigma_1' := \Cone (S,E_1) = 
\begin{pmatrix}
1 & 1 \\
1 & 0
\end{pmatrix}.
$$
On the toric chart $\bc^2_{\sigma_1'}$, the toric modification is written as
$\hat{\pi}_{\sigma_1'}(u_1,u_2) = (u_1u_2, u_1)$. 
The toric chart $\bc^2_{\sigma_1'}$ intersects the exceptional divisor $\hat{E}(S)$ only. 
If $\bbu_{\sigma_1'}^0 \in \tilde V \cap \hat{\pi}^{-1}(\zero) \cap \bc^2_{\sigma_1'}$, then 
$\bbu_{\sigma_1'}^0 \in \tilde V \cap \hat{E}(S)$. 
We have $f_S = z_2^2 \bar{z}_2$. 
Here we set 
$$r_1 := \rdeg_S f_S = 1\cdot 0 + 1\cdot (2+1) = 3
, \ \ \ \text{and} \ \ \ 
p_1 := \pdeg_S f_S = 1\cdot 0 + 1\cdot (2-1) = 1.
$$
Then we have 
$$
\frac{r_1 + p_1}{2} = \frac{3+1}{2} = 2
, \ \ \ \text{and} \ \ \ 
\frac{r_1 - p_1}{2} = \frac{3-1}{2} = 1, 
$$
and 
$$
\begin{array}{ccl}
(\hat{\pi}_{\sigma_1'}^* f)(u_1,u_2) & = & u_1^2\bar{u}_1 - 6u_1^3\bar{u}_1u_2^2 + 11u_1^3\bar{u}_1^2u_2^2\bar{u}_2^2 - 6u_1^4\bar{u}_1^2u_2^4\bar{u}_2^2 \\
                                               & = & u_1^2\bar{u}_1 (1-6u_1u_2^2 + 11u_1\bar{u}_1u_2^2\bar{u}_2^2 -6u_1^2\bar{u}_1u_2^4\bar{u}_2^2) .
\end{array}
$$
We set 
$$
\widetilde f(u_1,u_2) = 1-6u_1u_2^2 + 11u_1\bar{u}_1u_2^2\bar{u}_2^2 -6u_1^2\bar{u}_1u_2^4\bar{u}_2^2
$$
and 
the strict transform $\tilde V$ of $V$ to $X$ in the toric chart $\bc^2_{\sigma_1'}$ is given by
$$
\tilde V = \{ \widetilde f(u_1,u_2) = 0 \} .
$$
If $\bbu_{\sigma_1'}^0 \in \tilde V \cap \hat{\pi}^{-1}(\zero) \cap \bc^2_{\sigma_1'}$, then 
$\bbu_{\sigma_1'}^0 \in \tilde V \cap \hat{E}(S)$. 
Hence, $\bbu_{\sigma_1'}^0 = (0, u_2)$ for some $u_2 \in \bc$. 
However, we have $f(0, u_2) = 1 \neq 0$. Thus we have $\tilde V \cap \hat{E}(S) = \emptyset$ on the toric chart $\bc^2_{\sigma_1'}$. 
Thus, the Assumption (*) in Theorem 32 of \cite{Saito-Takashimizu2021winter} is satisfied for $\Sigma^*$ 
on the toric chart $\bc^2_{\sigma_1'}$. 

\medskip

\noindent
{\bf (II)}\ \ 
We next consider the cone 
$\sigma_2 := \Cone (S,P) = 
\begin{pmatrix}
1 & 1 \\
1 & 2
\end{pmatrix}
$ 
and the toric chart $\bc^2_{\sigma_2}$. 
Both $S$ and $P$ are strictly positive, and 
$\bc^2_{\sigma_2}$ intersects the exceptional divisors $\hat{E}(S)$ and $\hat{E}(P)$. 
If $\bbu_{\sigma_2}^0 \in \tilde V \cap \hat{\pi}^{-1}(\zero) \cap \bc^2_{\sigma_2}$, then 
$\bbu_{\sigma_2}^0 \in \tilde V \cap \hat{E}(S)$ or $\bbu_{\sigma_2}^0 \in \tilde V \cap \hat{E}(P)$. 
If $\bbu_{\sigma_2}^0 \in \hat{E}(S) \cap \hat{E}(P)$, then $\bbu_{\sigma_2}^0 = (0,0)$. 
If $\bbu_{\sigma_2}^0 \in \hat{E}(S)$ and $\bbu_{\sigma_2}^0 \not\in \hat{E}(P)$, then $\bbu_{\sigma_2}^0 = (0,u_2)$ for some $u_2 \ (\neq 0) \in \bc$. 
If $\bbu_{\sigma_2}^0 \not\in \hat{E}(S)$ and $\bbu_{\sigma_2}^0 \in \hat{E}(P)$, then $\bbu_{\sigma_2}^0 = (u_1,0)$ for some $u_1 \ (\neq 0) \in \bc$. 
Then $\bbu_{\sigma_2}^0$ is written as $(0,u_2)$ for some $u_2 \ (\neq 0) \in \bc$ 
in the toric chart $\bc^2_{\sigma_2'}$, where we set 
$$
\sigma_2' := \Cone (P,S) = 
\begin{pmatrix}
1 & 1 \\
2 & 1
\end{pmatrix}. 
$$
Thus, the Assumption (*) in Theorem 32 of \cite{Saito-Takashimizu2021winter} is satisfied for $\Sigma^*$ 
on the toric chart $\bc^2_{\sigma_2}$. 

\medskip

\noindent
{\bf (III)}\ \ 
We finally consider the cone 
$\sigma_3 := \Cone (P,E_2) = 
\begin{pmatrix}
1 & 0 \\
2 & 1
\end{pmatrix}
$ 
and the toric chart $\bc^2_{\sigma_3}$. 
On the toric chart $\bc^2_{\sigma_3}$, the toric modification is written as 
$\hat{\pi}_{\sigma_3}(u_1,u_2) = (u_1, u_1^2u_2)$. 
The toric chart $\bc^2_{\sigma_3}$ intersects the exceptional divisor $\hat{E}(P)$ only. 
If $\bbu_{\sigma_3}^0 \in \tilde V \cap \hat{\pi}^{-1}(\zero) \cap \bc^2_{\sigma_3}$, then 
$\bbu_{\sigma_3}^0 \in \tilde V \cap \hat{E}(P)$. 
Here $f_P = f$ and recall that 
$$r_1 := \rdeg_P f = 6
, \ \ \ \text{and} \ \ \ 
p_1 := \pdeg_S f_S = 2.
$$
Then we have 
$$
\frac{r_1 + p_1}{2} = \frac{6+2}{2} = 4
, \ \ \ \text{and} \ \ \ 
\frac{r_1 - p_1}{2} = \frac{6-2}{2} = 2, 
$$
and 
$$
\begin{array}{ccl}
(\hat{\pi}_{\sigma_3}^* f)(u_1,u_2) & = & u_1^4\bar{u}_1^2u_2^2\bar{u}_2 - 6u_1^4\bar{u}_1^2u_2\bar{u}_2 + 11u_1^4\bar{u}_1^2u_2 - 6u_1^4\bar{u}_1^2 \\
                                               & = & u_1^4\bar{u}_2 (u_2^2\bar{u}_2 - 6u_2\bar{u}_2 + 11u_2 - 6). 
\end{array}
$$
We set 
$$
\widetilde f(u_1,u_2) = u_2^2\bar{u}_2 - 6u_2\bar{u}_2 + 11u_2 - 6, 
$$
then 
the strict transform $\tilde V$ of $V$ to $X$ in the toric chart $\bc^2_{\sigma_3}$ is given by
$$
\tilde V = \{ \widetilde f(u_1,u_2) = 0 \} .
$$
If $\bbu_{\sigma_3}^0 \in \tilde V \cap \hat{\pi}^{-1}(\zero) \cap \bc^2_{\sigma_3}$, then 
$\bbu_{\sigma_3}^0 \in \tilde V \cap \hat{E}(P)$. 
Hence, $\bbu_{\sigma_3}^0 = (0, u_2)$ for some $u_2 \in \bc$. 
Since $\widetilde f(0,0)= -6 \neq 0$, we see that $u_2\neq 0$. 
Thus, the Assumption (*) in Theorem 32 of \cite{Saito-Takashimizu2021winter} is satisfied for $\Sigma^*$ 
on the toric chart $\bc^2_{\sigma_3}$. 

\bigskip

By the above arguments {\bf (I)}$\sim${\bf (III)}, 
the Assumption (*) in Theorem 32 of \cite{Saito-Takashimizu2021winter} 
is satisfied for the mixed polynomial germ $(f_{2,2,1,2,1,4}\ (k=3),\zero)$ and $\Sigma^*$. 
Hence, by Theorem 32, 
it is concluded that 
the strict transform $\tilde V$ of $V := f^{-1}(0)$ via the toric modification $\hat{\pi} : X \to \bc^2$ 
is a real analytic manifold outside of $\tilde V \cap \hat{\pi}^{-1}(\zero)$, and 
a topological manifold as a germ at $\tilde V \cap \hat{\pi}^{-1}(\zero)$. 

\bigskip

Moreover, we have the following theorem:
\begin{theorem}\label{real-analytic-manifold}
We set $f := f_{2,2,1,2,1,4}\ (k=3)$. 
The strict transform $\tilde V$ of $V := f^{-1}(0)$ via the toric modification $\hat{\pi} : X \to \bc^2$ 
is a real analytic manifold as a germ at $\tilde V \cap \hat{\pi}^{-1}(\zero)$. 
\end{theorem}

\begin{proof}
Recall the definition of $L(\Sigma^*)$ in Theorem 32 of \cite{Saito-Takashimizu2021winter}. 
It is sufficient to prove that $L(\Sigma^*) = \emptyset$. 
The cones of $\Sigma^*$ whose vertices are all strictly positive are 
$$\tau_1 := \Cone (S),\ \ \tau_2 := \Cone (P),\ \ \sigma_2 := \Cone (S,P).$$
For $\tau_1$, we have 
$\{ (\nu, \mu)\ |\ c_{\nu, \mu}\neq 0,\ \nu + \mu \not\in \De(S) \}  
= \{ ((2,1), (0,1)),\ ((2,1), (2,0)),\ ((4,0), (2,0)) \} 
\neq \emptyset$. 
Recall  that $\rdeg_S f_S = 3$. 
Hence, we have 
$$
\begin{array}{ccl}
\Lambda(\tau_1) & = & \min \{ S(\nu + \mu) - 3\ |\ (\nu, \mu) = ((2,1), (0,1)),\ ((2,1), (2,0)),\ ((4,0), (2,0)) \} \\
                        & =  & \min \{ {}^t(1,1)(\nu + \mu) - 3\ |\ \nu + \mu = (2,2), (4,1), (6,0) \} \\
                        & =  & \min \{ 4-3,\ 5-3,\ 6-3 \} = 1.
\end{array}
$$
For $\tau_2$, we have 
$\{ (\nu, \mu)\ |\ c_{\nu, \mu}\neq 0,\ \nu + \mu \not\in \De(P) \} = \emptyset$.  
For $\sigma_2$, we have 
$\{ (\nu, \mu, S)\ |\ c_{\nu, \mu}\neq 0,\ \nu + \mu \not\in \De(S) \} \cup \{ (\nu, \mu, P)\ |\ c_{\nu, \mu}\neq 0,\ \nu + \mu \not\in \De(P) \} 
= \{ ((2,1), (0,1), S),\ ((2,1), (2,0), S),\ ((4,0), (2,0), S) \} 
\neq \emptyset$. Hence, we have 
$\Lambda(\sigma_2) = \min \{ S(\nu + \mu) - 3\ |\ (\nu, \mu) = ((2,1), (0,1)),\ ((2,1), (2,0)),\ ((4,0), (2,0)) \} = 1$. 

On the other hand, we see that 
$\tilde V(\tau_1)^* \cap \bc_{\sigma_1}^2 = \emptyset,\ \tilde V(\tau_1)^* \cap \bc_{\sigma_2}^2 = \emptyset$, 
and moreover, 
$\tilde V(\sigma_2)^* \cap \bc_{\sigma_2}^2 = \emptyset$. 
Thus it is concluded that $L(\Sigma^*) = \emptyset$. 
This completes the proof of Theorem \ref{real-analytic-manifold}. 
\end{proof}

\bigskip

At the end of this paper, we present the following problem: 
\begin{problem}
\ \ 

\begin{itemize}
\item Does each $f_{2,2,1,2,1,4}\ (k\neq 3\ \text{and}\ k>2)$  have some mixed critical points on ${\bc^*}^2$?
\item More generally, 
does each $f_{2,b,c,d,e,4}$ in the cases II $\sim$ V ($a=2,\ \mathrm{f}=4$) of Table \ref{strongly-mixed-whp-5cases} 
have some mixed critical points on ${\bc^*}^2$?
\end{itemize}
\end{problem}
We need another useful criteria for mixed critical points like Proposition 1 of \cite{Oka2008}.



\begin{thebibliography}{99}

\bibitem{Gudkov74}
D. A. {\sc Gudkov}, 
The topology of real projective algebraic varieties, 
Usp. Mat. Nauk, {\bf 29--4} (1974), 3--79, 
English transl.,  
Russ. Math. Surveys, {\bf 29--4} (1974), 1--79. 

\bibitem{Ishi-Saito-Fukui2001}
G. {\sc Ishikawa}, S. {\sc Saito} {\sc and} T. {\sc Fukui} (the second part), I. {\sc Shimada} {\sc and} H. {\sc Tokunaga} (the first part), 
Algebraic curves and singularities (in Japanese), 
Kyoritsu shuppan, 2001. 


\bibitem{Kharlamov1976}
V. M. {\sc Kharlamov}, 
The topological type of nonsingular surfaces in $\br P^3$ of degree four, 
Funct. Anal. Appl. {\bf 10} (1976), 295--305. 

\bibitem{Nikulin79}
V. V. {\sc Nikulin}, 
Integral symmetric bilinear forms and some of their geometric applications, 
Izv. Akad. Nauk SSSR Ser Mat. 43-1 (1979), 111--177, 
English transl., 
Math. USSR Izv., {\bf 14}-1 (1980), 103--167. 

\bibitem{Oka1997}
M. {\sc Oka}, 
Non-degenerate complete intersection singularity, 
Hermann, Paris, 1997.

\bibitem{Oka2008}
M. {\sc Oka}, Topology of polar weighted homogeneous hypersurfaces, 
Kodai Math. J. {\bf 31} (2008), 163--182.

\bibitem{Oka2010}
M. {\sc Oka}, Non-degenerate mixed functions, 
Kodai Math. J. {\bf 33} (2010), 1--62.

\bibitem{Oka2015}
M. {\sc Oka}, Mixed functions of strongly polar weighted homogeneous face type, 
Singularities in Geometry and Topology 2011, 
Advanced Study in Pure Math. {\bf 66} (2015), 173--202.

\bibitem{Oka2018}
M. {\sc Oka}, 
Introduction to Complex and Mixed hypersurface singularities (in Japanese), 
Maruzen, Tokyo, 2018. 

\bibitem{Rokhlin72a}
V. A. {\sc Rokhlin}, 
Proof of Gudkov conjecture, 
Funkts. Anal. Prilozhen., {\bf 6--2} (1972), 62--64.  
English transl.,  
Funct. Anal. Appl., {\bf 6} (1972), 62--64. 

\bibitem{Rokhlin72b}
V. A. {\sc Rokhlin}, 
Congruences modulo 16 in Hilbert's sixteenth problem I, 
Funkts. Anal. Prilozhen., {\bf 6--4} (1972), 
English transl.,  
Funct. Anal. Appl., {\bf 6} (1972), 301--306. 

\bibitem{Rokhlin73}
V. A. {\sc Rokhlin}, 
Congruences modulo 16 in Hilbert's sixteenth problem II, 
Funkts. Anal. Prilozhen., {\bf 7--2} (1973), 91--92. 
English transl.,  
Funct. Anal. Appl., {\bf 7} (1973), 163--164. 

\bibitem{Saito-Takashimizu2021winter}
S. {\sc Saito} {\sc and} K. {\sc Takashimizu}, 
Resolutions of Newton non-degenerate mixed polynomials of 
strongly polar non-negative mixed weighted homogeneous face type, 
Kodai Math. J. {\bf 44} (2021), 457--491. 


\bibitem{Saito-Takashimizu2022}
S. {\sc Saito} {\sc and} K. {\sc Takashimizu}, 
A note on Newton non-degeneracy of mixed weighted homogeneous polynomials, 
arXiv:2107.08691, preprint. 


\bibitem{Viro1984}
O. Ya. {\sc Viro}, 
Gluing of plane algebraic curves and constructions of curves of degree 6 and 7, 
Lecture Notes in Math., {\bf 1060} (1984), 187--200. 

\bibitem{Viro1989}
O. Ya. {\sc Viro}, 
Real algebraic plane curves: Constructions with controlled topology, 
Algebra and analysis, {\bf 1--5} (1989); 
English transl. Leningrad Math. J., {\bf 1--5} (1990), 1059--1134.

\bibitem{Viro2006}
O. Ya. {\sc Viro}, 
Patchworking real algebraic varieties, arXiv:math/0611382 [math.AG], 2006. 

\bibitem{Wilson} 
G. {\sc Wilson}, 
Hilbert's sixteenth problem, 
Topology, {\bf 17} (1978), 
53--73. 

\end{thebibliography}
\end{document}